\documentclass[11pt,reqno]{amsart}
\usepackage{amssymb,amsfonts,mathrsfs,bm,enumerate,dsfont,mathptmx,amsaddr,color}

\makeatletter

\newtheorem{theorem}{Theorem}[section]
\newtheorem{lemma}[theorem]{Lemma}
\newtheorem{corollary}[theorem]{Corollary}
\newtheorem{proposition}[theorem]{Proposition}
\theoremstyle{definition}
 \newtheorem{remark}[theorem]{Remark}

\numberwithin{equation}{section}

\def\Ric{{\operatorname{Ric}}} 
 \def\II{{\operatorname{II}}}
 \def\Sect{{\operatorname{Sect}}}
\def\cut{\operatorname{cut}}

\newcommand\1{\hbox{\kern.375em\vrule height1.57ex depth-.1ex
    width.05em\kern-.375em \rm 1}}

 \newcommand\E{\mathbb{E}}
 
 \newcommand\R{\mathbb{R}}

  \newcommand\bull{{\hbox{\bf .}\thinspace}}
\newcommand\newdot{{\kern.8pt\cdot\kern.8pt}}

 \def\vd{\mathrm{d}} \def\r{\right}
\def\l{\left} \def\e{\operatorname{e}}

  \makeatother

\begin{document}

  \title[Functional inequalities on manifolds with non-convex boundary]
  {Functional inequalities on manifolds with non-convex
  boundary}

  \author[L.J.~Cheng, A.~Thalmaier and J.~Thompson]{Li-Juan Cheng\textsuperscript{1,2}, Anton Thalmaier\textsuperscript{1} and James Thompson\textsuperscript{1}}

  \address{\textsuperscript{1}Mathematics Research Unit, FSTC, University of Luxembourg\\
    Maison du Nombre, 6, Avenue de la Fonte,\\ 4364 Esch-sur-Alzette, Grand Duchy of Luxembourg}
  \address{\textsuperscript{2}Department of Applied Mathematics, Zhejiang University of Technology\\
    Hangzhou 310023, The People's Republic of China}
  \email{lijuan.cheng@uni.lu \text{and} chenglj@zjut.edu.cn}
  \email{anton.thalmaier@uni.lu}

\begin{abstract}
In this article, new curvature conditions are introduced to establish functional inequalities including gradient estimates, Harnack inequalities and transportation-cost inequalities on manifolds with non-convex boundary.
\end{abstract}

\date{\today}
\maketitle

\section{Introduction}

Let $(M,g)$ be a $d$-dimensional complete and connected Riemannian manifold with Riemannian distance $\rho$, boundary $\partial M$ and inward pointing unit normal vector $N$. Define the second fundamental form of the boundary by
$$\II(X,Y)=-\l<\nabla_XN,Y\r>, \quad X,Y\in T_x\partial M,\quad x \in \partial M$$
where $T\partial M$ denotes the tangent bundle of $\partial M$. In order to study non-convex boundaries, we will perform a conformal change of metric such that the boundary is convex under the new metric. In particular, we will use the fact that if
$$\mathscr{D}:=\{\phi\in C_b^2(M)\colon \inf \phi = 1, \II \geq - N\log \phi \}$$
and $\phi\in \mathscr{D}$ then the boundary $\partial M$ is convex under the metric $\phi^{-2}g$ (see \cite[Theorem 1.2.5]{Wbook2}).

Given a $C^1$-vector field $Z$ on $M$, consider the elliptic operator $L:=\Delta+Z$ and let $X^x_t$ be a reflecting $L$-diffusion process starting from $X_0^x=x$. Then $X_t^x$ solves the Stratonovich equation
\begin{align*}
\vd X_t^x=\sqrt{2} u_t^x\circ \vd B_t+Z(X_t^x)\,\vd t+N(X_t^x)\,\vd l_t^x, \quad X^x_0 =x
\end{align*}
where $u_t^x$ is the horizontal lift of $X_t^x$ to the orthonormal frame bundle ${\rm O}(M)$ with $\pi(u_0^x)=x$, $B_t$ is a standard $\R^d$-valued Brownian motion defined on a complete naturally filtered probability space $(\Omega, \{\mathscr{F}_t\}_{t\geq 0},\mathbb{P})$ and $l_t^x$ is a continuous adapted nondecreasing and nonnegative process which increases only on $\lbrace t \geq 0\colon X_t^x \in \partial M \rbrace$. The process $l_t^x$ is the local time of $X^x_t$ on $\partial M$.

We assume that $X_t^x$ is non-explosive for each $x\in M$. Then the diffusion process $X_t^x$ gives rise to the Neumann semigroup $P_t$ which solves the diffusion equation $(\partial_t -L)P_t = 0$ with Neumann boundary condition $NP_t = 0$. Furthermore $P_tf(x)=\E[f(X_t^x)]$ for each $f\in C_b(M)$ .

In \cite{Hsu}, Hsu found a probabilistic formula for $\nabla P_tf$ for compact manifolds with boundary, which he used to derive a gradient estimate. Feng-Yu Wang extended it to the non-compact case \cite[Theorem 3.2.1]{Wbook2} under the assumption that $|\nabla P_{\bull}f|$ is uniformly bounded on $[0,t]\times M$. Wang's formula is given below by Theorem \ref{Bismut-formula}. In \cite[Proposition 3.2.7]{Wbook2}, he proved that if $\Ric^Z:=\Ric-\nabla Z\geq K$ for some $K\in C(M)$ and there exists $\phi\in \mathcal{D}$ such that
\begin{equation}\label{eq:wangconst}
\tilde{K}_{\phi}:=\inf_M\l\{\phi^2 K+\frac{1}{2}L\phi^2-|\nabla \phi^2|\,|Z|-(d-2)|\nabla \phi|^2\r\}>-\infty
\end{equation}
then $|\nabla P_{\bull}f|$ is uniformly bounded on $[0,t]\times M$.

In this article, we revisit this problem by using coupling methods to weaken this curvature condition. In particular, we prove (see Theorem \ref{main-th1}) that if there exists $\phi\in \mathcal{D}$ and a constant $K_{\phi}$ such that
$$\Ric^Z+L\log \phi-|\nabla \log \phi|^2\geq K_{\phi}$$
then $|\nabla P_{\bull}f|$ is uniformly bounded on $[0,t]\times M$. The upper bound we obtain improves that of \cite[Proposition 3.2.7]{Wbook2}. We construct such a function $\phi$ in Proposition \ref{sec3-prop}, under the assumption that there exist non-negative constants $\sigma$ and $\theta$ such that $-\sigma\leq \II \leq \theta$ and a positive constant $r_0$ such that on $\partial_{r_0}M:=\{x\in M\colon \rho_{\partial}(x)\leq r_0\}$ the function $\rho_{\partial}$ is smooth, the norm of $Z$ is bounded and $\Sect \leq k$ for some positive constant $k$.

F.-Y. Wang also used coupling methods to consider Harnack and transportation-cost inequalities on manifolds with boundary \cite{Wbook2}. We reconsider these problems too and find that the curvature conditions used to establish these inequalities can also be weakened and simplified. It is worth mentioning that we find a transportation-cost inequality on the path space of the reflecting diffusion process which (see Theorem \ref{thm-tran-ineq}) recovers the results for the convex boundary case, making the theory of functional inequalities on path space complete.

This article is organized as follows. In Section 2, we prove the gradient estimates, Harnack inequalities and transportation-cost inequalities for the Neumann semigroup via coupling methods. In Section 3, we construct a function $\phi$ which satisfies the new curvature conditions.

\section{Functional inequalities}

\subsection{Gradient estimates}

The derivative formula for $P_t$ is known as Bismut formula (see \cite{Bismut,EL}). The formula we introduce is of a more general type due to Thalmaier \cite{Thalmaier}. As mentioned in the introduction, Hsu \cite{Hsu} found this type of formula for compact manifolds with boundary. The following formula for manifolds with boundary, due to F.-Y. Wang \cite[Theorem~3.2.1]{Wbook2}, does not require compactness. See also \cite{Arnaudon_Li:2017} for recent work on
probabilistic representations of the derivative of Neumann semigroups.

\begin{theorem}\label{Bismut-formula}
Let $t>0$ and $u_0\in {\rm O}_x(M)$ be fixed. Let $K\in C(M)$ and $\sigma\in C(\partial M)$ be such that $\Ric_Z \geq K$ and $\II \geq \sigma$. Assume that
\begin{align*}
\sup_{s\in [0,t]} \E^x \left[\exp {\l(-\int_0^sK(X_r)\,\vd r-\int_0^s\sigma(X_r)\,\vd l_r\r)}\right]<\infty.
\end{align*}
Then there exists a progressively measurable process $\{Q_s\}_{s\in [0,t]}$ on $\R^d\otimes\R^d$
such that
\begin{align*}
Q_0=I,\quad \|Q_s\|\leq \exp\l(-\int_0^sK(X_r)\,\vd r-\int_0^s\sigma(X_r)\,\vd l_r\r),\quad s\in [0,t]
\end{align*}
and for any $f\in C_b^1(M)$ such that $\nabla P_{\bull}f$ is bounded on $[0,t]\times M$, for any $h\in C^1([0,t])$ with $h(0)=0$ and $h(t)=1$, we have
\begin{align*}
u_0^{-1}\nabla P_tf(x)&=\E^x\left[Q_tu_t^{-1}\nabla f(X_t)\right]=\frac1{\sqrt{2}}\,\E^x\left[f(X_t)\int_0^t \dot{h}(s)Q_s\,\vd B_s\right].
\end{align*}
\end{theorem}

In order to use this formula it is necessary to check the uniform boundedness of $\nabla P_{\bull}f$ on $[0,t]\times M$. In \cite[Proposition 3.2.7]{Wbook2}, F.-Y. Wang did so using a conformal change of metric such that under the new metric the boundary is convex, and by then making a time change of the $L$-diffusion process $X_t$. Here, we use coupling methods to study this problem again and obtain improved upper bounds.

\begin{theorem}\label{main-th1}
If there exist $\phi\in \mathscr{D}$ and a constant $K_{\phi}$ such that
\begin{align}\label{K-phi}
\Ric^Z+L \log \phi- |\nabla \log \phi|^2 \geq K_{\phi}
\end{align}
then for $f\in C^1(M)$ such that $f$ is constant outside a compact set,
$$|\nabla P_{t}f|\leq \|\phi\|_{\infty}\|\nabla f\|_{\infty} \e^{-K_{\phi}t},\quad t>0.$$
\end{theorem}

\begin{proof}
We use a coupling method to prove the gradient inequality. To this end, we need to conformally change the metric $g$. Since $\phi \in \mathscr{D}$, the boundary $\partial M$ is convex under new metric $g':=\phi^{-2}g$. Let ${\Delta}'$ and ${\nabla}'$ be the Laplacian and gradient operator associated with the metric $g'$. Then
\begin{align}\label{eq-add-3}
L&=\phi^{-2}\l({\Delta}'+\phi^2\l(Z+(d-2)\nabla \log \phi\r)\r)\notag\\
&=\phi^{-2}\l(\sum_{i=1}^d V_i^2+\phi^2\l(Z+(d-2)\nabla \log \phi\r)\r)\notag\\
&=\sum_{i=1}^d(\phi^{-1}V_i)^2+\phi^{-2}Z',
\end{align}
where $\{V_i\}_{i=1}^d$ is a $g'$-orthonormal basis of $T_xM$ and
$Z':=\phi^2\left(Z+(d-1)\nabla \log \phi\right)$ (the factor $d-1$ corrects the factor $d-2$ appearing in the proof of \cite[Proposition 3.2.7]{Wbook2}).
Now consider the process $X_t$ generated by $L$ on the manifold $(M,g')$. Denoting by $\vd_I$ the It\^{o} differential on $M$, it is easy to see that $X_t$ is a solution to the equation
\begin{align}\label{add-eq-2}
\vd_IX_t=\sqrt{2}\phi^{-1}(X_t)u_t\,\vd B_t+\phi^{-2}(X_t)Z'(X_t)\,\vd t+N'(X_t)\,\vd l_t,\quad X_{0}=x
\end{align}
where the horizontal lift $u_t$%, Brownian motion $B_t$
and boundary local time $l_t$ are now defined with respect to the metric $g'$. Recall that in local coordinates, the It\^{o} differential of a continuous semimartingale $X_t$ on $M$ is given (see \cite{Emery} or \cite{ATW2006}) by
\begin{align*}
(\vd _I X_t)^k=\vd X_t^k+\frac{1}{2}\sum_{i,j=1}^d{\Gamma'}_{ij}^k(X_t)\,\vd \langle X^i, X^j \rangle_t,\quad 1\leq k\leq d
\end{align*}
where ${\Gamma'}_{ij}^k$ are the Christoffel symbols of $g'$. Similarly, let $Y_t$ solve
\begin{align}\label{add-eq-3}
\vd_IY_t=\sqrt{2}\phi^{-1}(Y_t)\tilde{u}_t\,\vd \tilde{B}_t+\phi^{-2}(Y_t)Z'(Y_t)\,\vd t+N'(Y_t)\,\vd \tilde{l}_t,\quad Y_{0}=y
\end{align}
with horizontal process $\tilde{u}_t$ and boundary local time $\tilde{l}_t$ and where $\tilde{B}_t$ is a new Brownian motion satisfying
$$1_{\{(X_t,Y_t)\notin\cut\}}\tilde{u}_t\,\vd \tilde{B}_t=1_{\{(X_t,Y_t)\notin\cut\}}P'_{X_t,Y_t}u_t\,\vd B_t$$
with $\cut\subset M\times M$ denoting the set of cut points.
For the sake of conciseness, we may assume without loss of generality that the cut locus of $M$ is empty. Denote by $\rho'$ the distance function for the metric~$g'$. Since the boundary $\partial M$ is convex under~$g'$, by the It\^{o} formula, we have
\begin{align*}
\vd \rho'(X_t,Y_t)&\leq \sqrt{2}\left(\phi^{-1}(X_t)-\phi^{-1}(Y_t)\right)\vd b_t+
\bigg\{\sum_{i=1}^{d}(U_i)^2\rho'(X_t,Y_t)\\
&\quad+\l<\phi^{-2}(Y_t)Z'(Y_t),{\nabla}'\rho'(X_t,\newdot)(Y_t)\r>'+\l<\phi^{-2}(X_t)Z'(X_t),{\nabla}'\rho'(\newdot, Y_t)(X_t)\r>'\bigg\}\,\vd t
\end{align*}
where $b_t$ is a one-dimensional Brownian motion, $\{U_i\}_{i=1}^{d}$ are vector fields on $M\times M$ such
that ${\nabla}' U_i(X_t,Y_t)=0$ and
$$U_i(X_t,Y_t)=(\phi^{-1}(X_t)V_i,\phi^{-1}(Y_t)\,P_{X_t,Y_t}'V_i),\quad1\leq i\leq d$$
for $\{V_i\}_{i=1}^{d}$ a $g'$-orthonormal basis of $T_{X_t}M$. Here $P_{x,y}'$ denotes parallel displacement from $x$ to $y$ with respect to the metric $g'$. Write ${\rho'}={\rho'}(X_t,Y_t)$ and for a minimizing $g'$-geodesic $\gamma$ with $\gamma(0) = X_t$ and $\gamma(\rho') =Y_t$ let
$${J_i}(s)=\phi^{-1}(\gamma(s))\,{P}'_{\gamma(0),\gamma(s)}V_i,\quad1\leq i\leq d$$
where ${J_i}(0)=\phi^{-1}(X_t)V_i$ and ${J_i}({\rho'})=\phi^{-1}(Y_t){P}'_{X_t,Y_t}V_i$. Since ${P}'_{\gamma(0),\gamma(s)}V_i$ are parallel vector fields along $\gamma $ with respect to the metric $g'$, we have
\begin{align}\label{add-eq-7}
&\sum_{i=1}^{d}(U_i)^2\rho'(X_t,Y_t)\notag\\
&\leq \sum_{i=1}^{d}\int_0^{\rho'}\l\{|{\nabla}'_{\dot{\gamma}}J_i|'^2-\l<R'(\dot{\gamma},J_i)J_i, \dot{\gamma}\r>'\r\}(s)\,\vd s\nonumber\\
&= d\int_0^{\rho'}\phi^{-2}(\gamma(s))\langle \nabla \log \phi(\gamma(s)),\dot{\gamma}(s)\rangle^2\,\vd s-\int_0^{\rho'}\phi^{-2}(\gamma(s))\Ric'(\dot{\gamma}(s),\dot{\gamma}(s))\,\vd s.
\end{align}
On the other hand
\begin{align}\label{add-eq-4}
&\phi^{-2}(X_t)\l<Z'(X_t), {\nabla}'\rho'(\newdot, Y_t)(X_t)\r>'+\phi^{-2}(Y_t)\l<Z'(Y_t), {\nabla}'\rho'(X_t,\newdot)(Y_t)\r>'\nonumber\\
&=
\int_0^{\rho'}\frac{\vd }{\vd s}\l\{\phi^{-2}(\gamma(s))\l<Z'({\gamma(s)}),\dot{\gamma}(s)\r>'\r\}\,\vd s\nonumber\\
&= \int_0^{\rho'}\phi^{-2}(\gamma(s))\l<({\nabla}'_{\dot{\gamma}}Z')\circ \gamma, \dot{\gamma}\r>'(s)\,\vd s\notag\\
&\quad-2\int_0^{\rho'}\phi^{-2}(\gamma(s))\langle \nabla \log \phi(\gamma(s)),\dot{\gamma}(s)\rangle\l<Z'({\gamma(s)}),\dot{\gamma}(s)\r>'\,\vd s.
\end{align}
Moreover
\begin{align*}
\langle Z'(\gamma(s)),\dot{\gamma}(s) \rangle'=\langle Z,\dot{\gamma}(s)\rangle+(d-1)\langle\nabla \log \phi, \dot{\gamma}(s)\rangle.
\end{align*}
Combining this with \eqref{add-eq-7} and \eqref{add-eq-4}, we have
\begin{align}\label{add-eq-5}
\vd \rho'&(X_t,Y_t)\leq \sqrt{2}(\phi^{-1}(X_t)-\phi^{-1}(Y_t))\,\vd b_t\notag\\
&
-\int_0^{\rho'}\phi^{-2}[(\Ric^Z)'(\dot{\gamma}(s),\dot{\gamma}(s))+(d-3)\langle\nabla \log \phi, \dot{\gamma}(s)\rangle^2+2\langle \nabla \log \phi, \dot{\gamma}(s)\rangle \langle Z,\dot{\gamma}(s) \rangle]\,\vd s.
\end{align}
By \cite[Theorem~1.159]{Besse} we know that
\begin{align*}
(\Ric^{Z})'(\dot{\gamma},\dot{\gamma})&=\Ric'(\dot{\gamma},\dot{\gamma})-\langle{\nabla}_{\dot{\gamma}}{Z}',\dot{\gamma}\rangle'\\
&=\Ric^Z(\dot{\gamma},\dot{\gamma})+\frac{1}{2}L\phi^2-2\l<\nabla\log\phi, \dot{\gamma}\r>\l<Z,\dot{\gamma}\r> -(d-2)\l<\dot{\gamma},\nabla\log\phi\r>^2 - 2|\nabla \phi|^2
\end{align*}
and, noting that $|\dot{\gamma}|=\phi$, we thus have
\begin{align}\label{Ric}
(\Ric^Z)'&(\dot{\gamma}(s),\dot{\gamma}(s))+(d-3)\langle\nabla \log \phi, \dot{\gamma}(s)\rangle^2+2\langle \nabla \log \phi, \dot{\gamma}(s)\rangle \langle Z,\dot{\gamma}(s) \rangle\notag\\
=\text{ }&\Ric^Z(\dot{\gamma}(s),\dot{\gamma}(s))+\frac{1}{2}L\phi^2-\l<\dot{\gamma},\nabla\log\phi\r>^2\notag - 2|\nabla \phi|^2 \\
\geq\text{ } & \Ric^Z(\dot{\gamma}(s),\dot{\gamma}(s))+\frac{1}{2}L\phi^2-3|\nabla \phi|^2\notag\\
=\text{ }& \Ric^Z(\dot{\gamma}(s),\dot{\gamma}(s))+\phi^2 L\log \phi-|\nabla \phi|^2.
\end{align}
Consequently, if
\begin{align*}
\Ric^Z+(L\log \phi-|\nabla \log \phi|^2) \langle \cdot,\cdot \rangle \geq K_{\phi}\langle \cdot,\cdot \rangle
\end{align*}
then, combining \eqref{add-eq-5} with \eqref{Ric}, we arrive at
\begin{align*}
\vd {\rho}'(X_t,Y_t)\leq \sqrt{2}(\phi^{-1}(X_t)-\phi^{-1}(Y_t)) \,\vd b_t-K_{\phi}{\rho'(X_t, Y_t)}\,\vd t.
\end{align*}
Then, observing that $\rho'\leq \rho\leq \|\phi\|_{\infty}\rho'$, we have
\begin{align*}
|\nabla P_{t}f|(x)&=\lim_{y\rightarrow x}\bigg|\frac{P_{t}f(x)-P_{t}f(y)}{\rho(x,y)}\bigg|\\
&=\lim_{y\rightarrow x}\bigg|\E^{(x,y)}\l[\frac{f(X_{t})-f(Y_{t})}{\rho(X_{t},Y_{t})}\frac{\rho(X_{t},Y_{t})}{\rho'(X_{t},Y_{t})}\frac{\rho'(X_{t},Y_{t})}{\rho'(x,y)}\frac{\rho'(x,y)}{\rho(x,y)}\r]\bigg|\\
&\leq \|\phi\|_{\infty}\|\nabla f\|_{\infty}\e^{-K_{\phi}t}
\end{align*}
which completes the proof.
\end{proof}

\begin{remark}\label{rem-1}\
\begin{enumerate}[(i)]
\item Since $(U_d)^2\rho'\neq 0$, it was indeed necessary to account for this quantity in inequality \eqref{add-eq-7}, correcting the proof of \cite[Theorem 3.4.6]{Wbook2}.
\item Compared with the proof of \cite[Theorem 3.4.6]{Wbook2}, our choice of vector field $J_i$ yields a simpler result.
\item In \cite{Wa14}, a certain technical assumption which was used to ensure the uniformly boundedness of $|\nabla P_{\bull}f|$ on $[0,t]\times M$ is no longer needed in the results.
\end{enumerate}
\end{remark}

The following results remove the additional condition in \cite[Corollary 3.6.5 (1)]{Wbook2} or \cite[Corollary 1.2 (1)]{Wa14} to ensure the uniform boundedness of $|\nabla P_{\bull}f|$ on $[0,t]\times M$ and give a another proof to extend these inequalities to $L^p$ forms for $p>1$:

\begin{theorem}\label{gradient-main-2}
If there exists $\phi\in \mathscr{D}$ such that for $p> 1$ the inequality
\begin{align}\label{eq:bndone}
\Ric^Z+L\log \phi-p|\nabla \log \phi|^2\geq K_{\phi,p}
\end{align}
holds, then for $t>0$ and $f\in C^1_b(M)$,
$$|\nabla P_{t}f|\leq \frac{1}{\phi}\e^{-K_{\phi,p}t}\left(P_{t}(\phi|\nabla f|)^{p/(p-1)}\right)^{(p-1)/p}.$$
\end{theorem}

\begin{proof}
The lower bound \eqref{eq:bndone} implies $\Ric^Z+L\log \phi-|\nabla \log \phi|^2$ is bounded below. By Theorem \ref{main-th1}, it follows that
$|\nabla P_{\bull}f|$ is bounded on $[0,t]\times M$. Furthermore
$$\Ric^Z\geq K_{\phi,p}-L\log \phi +p|\nabla \log \phi|^2=K_{\phi,p}+\frac{1}{p}\phi^{p}L\phi^{-p}\  \ \mbox{and}\quad\II \geq -N\log \phi$$
and so, by Theorem \ref{Bismut-formula}, there exists $\lbrace Q_s \rbrace_{s\in [0,t]}$ such that
\begin{align}\label{Q-upper}
\|Q_t\|\leq \exp\l(-K_{\phi,p} t-\frac{1}{p}\int_0^t\phi^{p}L\phi^{-p}(X_s)\,\vd s+\int_0^tN\log \phi(X_s)\,\vd l_s\r)
\end{align}
with
\begin{align}\label{gradient-p}
|\nabla P_tf|^p\leq \left(P_t(\phi |\nabla f|)^{p/(p-1)}\right)^{p-1}\E\left[ \phi^{-p}(X_t)\|Q_t\|^p\right].
\end{align}
It therefore suffices to give the upper bound estimate of the following term:
\begin{align*}
\E\l[\phi^{-p}(X_t)\exp\l(-\int_0^t\phi^{p}L\phi^{-p}(X_s)\,\vd s+p\int_0^tN\log \phi(X_s)\,\vd l_s\r)\r].
\end{align*}
To this end, by the It\^{o} formula, it is easy to see that
\begin{align*}
\vd \phi^{-p}(X_t)&=\langle \nabla \phi^{-p}(X_t), u_t\,\vd B_t \rangle+L\phi^{-p}(X_t)\,\vd t +N\phi^{-p}(X_t)\,\vd l_t\\
&=\langle \nabla \phi^{-p}(X_t), u_t\,\vd B_t \rangle-p\phi^{-p}(X_t)\l(-\frac{1}{p}\phi^{p}L\phi^{-p}(X_t)\,\vd t +N\log\phi(X_t)\,\vd l_t\r).
\end{align*}
So
\begin{align*}
M_t=\phi^{-p}(X_t)\exp\l(-\int_0^t\phi^{p}(X_s)L\phi^{-p}(X_s)\,\vd s+p\int_0^tN\log \phi(X_s)\,\vd l_s\r)
\end{align*}
is a local martingale. Thus,
\begin{align*}
\E\l[\phi^{-p}(X_t)\exp\left(-\int_0^t\phi^{p}(X_s)L\phi^{-p}(X_s)\,\vd s+p\int_0^tN\log \phi(X_s)\,\vd l_s\right)\r]\leq \phi^{-p}(x).
\end{align*}
Combining this and \eqref{gradient-p} and \eqref{Q-upper} completes the proof.
\end{proof}

\begin{corollary}\label{cor-1}
If there exists $\phi\in \mathscr{D}$ such that for $p> 1$ the inequality
\begin{align*}
\Ric^Z+L\log \phi-p|\nabla \log \phi|^2\geq K_{\phi,p}
\end{align*}
holds, then for $t>0$ and $f\in C^1_b(M)$,
$$|\nabla P_{t}f|\leq \|\phi\|_{\infty}\e^{-K_{\phi,p}t}(P_{t}|\nabla f|^{p/(p-1)})^{(p-1)/p};$$
 and for $f\in \mathscr{B}_b(M)$ and $t>0$,
\begin{align}\label{equ1}
|\nabla P_{t}f|^2\leq \|\phi\|_{\infty}^{2}\,\frac{K_{\phi,2}}{\e^{2K_{\phi,2}t}-1}P_{t}f^2.
\end{align}
\end{corollary}
\begin{proof}
The first assertion follows from Theorem \ref{gradient-main-2} by observing $\phi\geq 1$. As Theorem \ref{Bismut-formula} can be used under our condition directly,  the main idea of the proof of \eqref{equ1} is similar to that of \cite[Corollary 3.2.8]{Wbook2}, we skip it here.
\end{proof}

Note that taking the limit $p\downarrow 1$ in Corollary \ref{cor-1} recovers Theorem \ref{main-th1}.

\subsection{Harnack inequalities}
In \cite[Theorem 3.4.7]{Wbook2} or \cite[Theorem 3.1]{W11}, F.-Y. Wang used a coupling method to obtain dimension free Harnack inequalities and a log-Harnack inequality on manifolds with boundary assuming $\Ric^Z\geq K$ for some $K\in C(M)$ with $\phi\in \mathscr{D}$ such that $\tilde{K}_{\phi}$ is finite (where the quantity $\tilde{K}_{\phi}$ is defined as in \eqref{eq:wangconst}). The coefficient involved in these inequalities is: $$\tilde{\kappa}_{\phi}=2{{\tilde{K}}_{\phi}}^- +4\|\phi Z+(d-2)\nabla \phi\|_{\infty}\|\nabla \log \phi\|_{\infty}+2d\|\nabla \log \phi\|_{\infty}^2.$$ However, in \cite[Corollary 3.6.5 (2)]{Wbook2} or \cite[Corollary 1.2 (2)]{Wa14}, F.-Y. Wang used a modified curvature condition to obtain a log-Harnack inequality with coefficient $K_{\phi,2}$. A natural question raised is how to get  dimension-free Harnack inequalities with coefficient $K_{\phi,2}$ which are consistent with the log-Harnack inequality \cite[Corollary 3.6.5 (2)]{Wbook2} or \cite[Corollary 1.2 (2)]{Wa14}. We give the answer as follows, also by using coupling methods:

\begin{theorem}\label{Harnack-th2}
Assume there exists $\phi\in \mathcal{D}$ such that
\begin{align}\label{K-phi-1}
\Ric^Z+ L\log\phi-2|\nabla \log \phi|^2\geq {K}_{\phi,2}
\end{align}
for some  constant ${K}_{\phi,2}$. Then
\begin{enumerate}[\rm(i)]
\item for $T>0$, $x,y\in M$ and a measurable function $f\geq 1$, we have
\begin{align*}
P_T\log f(y)\leq \log P_Tf(x)+\frac{K_{\phi,2}\,\|\phi\|^2_{\infty}\,\rho^2(x,y)}{2(\e^{2K_{\phi,2}T}-1)};
\end{align*}
\item for $T>0$, $x,y\in M$, $p>\|\phi\|_{\infty}^2$ and function $f\in C^{1}_b(M)$, we have
$$\big(P_{T}f(y)\big)^p\leq P_{T}f^p(x)\exp\left(\frac{\sqrt{p}(\sqrt{p}-1)\,K_{\phi,2}\,\|\phi\|_{\infty}^2\,\rho^2(x,y)}{8\delta_p(\sqrt{p}-1-\delta_p)(\e^{2K_{\phi,2}T}-1)}\right),$$
where $\delta_{p}=\max\l\{\|\phi\|_{\infty}-1,\frac{\sqrt{p}-1}{2}\r\}.$
\end{enumerate}
\end{theorem}

\begin{proof}
Fix $x,y \in M$ and $T>0$. As in the proof of Theorem \ref{main-th1}, we consider the process $X_t$ generated by $L=\Delta+Z$ under the metric $g':=\phi^{-2}g$, for which the boundary $(\partial M, g')$ is convex. Let $X_t$ solve equation (\ref{add-eq-2}) with $X_0=x$. For a strictly positive function $\xi\in C([0,T))$, to be later determined, let $Y_t$ solve
\begin{align*}%\label{SDE2}
\vd _I Y_t=\sqrt{2}\phi^{-1}(Y_t)P'_{X_t, Y_t}u_t\,\vd B_t+\phi^{-2}_tZ'(Y_t)\,\vd
t-\frac{\phi^{-1}(Y_t)\rho'(X_t,Y_t)}{\phi^{-1}(X_t)\xi
(t)}\nabla'\rho'(X_t,\newdot)(Y_t)\,\vd
t+N'(Y_t)\,\vd \tilde{l}_t
\end{align*}
for $t\in [0,T)$ with $Y_0=y$, where $\tilde{l}_t$ is the local time of $Y_t$ on $\partial M$. Now consider the process $(X_t,Y_t)$ starting from $(x,y)$, which is a well defined continuous process for $t\leq T\wedge \zeta$ where $\zeta$ is the explosion time of $Y_t$; that is $\zeta:=\lim_{n\rightarrow \infty} \zeta_n$ for $\zeta_n:=\inf\{t>0: \rho'(y,Y_t)\geq n\}$. As in the proof of Theorem \ref{main-th1}, we can assume that the cut-locus of $(M,g')$ is empty, so that parallel displacement is smooth. Let
\begin{align}\label{e8}
\vd \tilde{B}_t=\vd
B_t+\frac{\rho'(X_t,Y_t)}{\sqrt{2}\xi(t)\phi^{-1}(X_t)}u_t^{-1}\nabla'\rho'(\newdot,Y_t)(X_t)\,\vd
t,\quad 0\leq t< T\wedge \zeta.
\end{align}
Since
\begin{align*}
\Ric^Z+L\log \phi-|\nabla \log \phi|^2\geq {K}_{\phi,2}+|\nabla \log \phi|^2,\quad t\in [0,T],
\end{align*}
by a similar calculation as for \eqref{add-eq-5} we find
\begin{align}\label{2eq6}
\vd \rho'(X_t,Y_t)&\leq
\sqrt{2}(\phi^{-1}(X_t)-\phi^{-1}(Y_t))\l<\nabla'\rho'(\newdot,Y_t)(X_t),u_t\vd
\tilde{B}_t\r>'\notag\\
&\quad-\l(\int_0^{\rho'(X_t,Y_t)}(K_{\phi,2}+|\nabla \log \phi|^2)(\gamma(s))\,\vd s\r)\vd t
-\frac{\rho'(X_t,Y_t)}{\xi(t)}\,\vd t,\quad 0\leq t<T\wedge \zeta
\end{align}
which implies
\begin{align}\label{rho-1}
\vd\frac{\rho'(X_t,Y_t)^2}{\xi(t)}&\leq
\frac{2\sqrt{2}}{\xi(t)}\rho'(X_t,Y_t)\l(\phi^{-1}(X_t)-\phi^{-1}(Y_t)\r)\l<\nabla'\rho'(\newdot,
Y_t)(X_t),u_t\vd
\tilde{B}_t\r>'\nonumber\\
&\quad-\frac{\rho'(X_t,Y_t)^2}{\xi^2(t)}\l(\dot{\xi}(t)+2K_{\phi,2}\xi(t)+2\r)\vd
t,\quad 0\leq t< T.
\end{align}
Now for $\theta\in (0,2)$ let
$$\xi(t)=(2-\theta)\int_t^{T}\e^{-2K_{\phi,2}(t-s)}\,\vd s,\quad t\in [0,T)$$
so that $\xi$ solves the equation
$$\dot{\xi}(t)+2K_{\phi,2}\xi(t)+2=\theta,\quad t\in [0,T).$$
Combining this with \eqref{rho-1}, we find
\begin{align*}%\label{rho-estimate}
\vd\frac{\rho'(X_t,Y_t)^2}{\xi(t)}\leq
\frac{2\sqrt{2}}{\xi(t)}\rho'(X_t,Y_t)\l(\phi^{-1}(X_t)-\phi^{-1}(Y_t)\r)\l<\nabla'\rho'(\newdot,
Y_t)(X_t),u_t\vd
\tilde{B}_t\r>'-\frac{\rho'(X_t,Y_t)^2}{\xi(t)^2}\theta\,\vd t.
\end{align*}
The remainder of argument is given by the proof of \cite[Theorem 3.4.7]{Wbook2}.
\end{proof}

\subsection{Transportation-cost inequalities}
Consider $\mu,\nu\in\mathscr{P}(M)$ where $\mathscr{P}(M)$ denotes the space of all probability measures on $M$. Recall the $L^p$-Wasserstein distance between $\mu$ and $\nu$ is
$$W_{p}(\mu,\nu)=\inf_{\eta\in \mathscr{C}(\mu,\nu)}\l\{\int_{M\times M}\rho(x,y)^p\,\vd \eta(x,y)\r\}^{1/p}$$
where $\mathscr{C}(\mu,\nu)$ is the set for couplings of $\mu$ and $\nu$. When the manifold has no boundary,
it is well known that the curvature condition,
$$\Ric^Z\geq K\quad \mbox{for some constant }K$$
is equivalent to
$$W_{p}(\mu P_t, \nu P_t)\leq W_p(\mu, \nu)\e^{-Kt},\quad \mu,\nu\in \mathscr{P}(M),$$
where $\mu P_t\in \mathscr{P}(M)$ is defined by $(\mu P_t)(A)=\mu(P_t1_A)$ for measurable set $A$.
 This equivalence is due to \cite{vonS05} which is extended to the manifolds with convex boundary \cite{W13}.
Using a coupling method, we obtain the following transportation-cost inequality.
\begin{theorem}\label{s3-t1}
If there exists $\phi\in \mathscr{D}$ and a constant $K_{\phi,2}$
satisfying
\begin{align*}
\Ric^Z+ L\log \phi-2|\nabla \log \phi|^2\geq K_{\phi,2}
\end{align*}
then
$$W_2(\mu P_t,\nu P_t)\leq \|\phi\|_{\infty}\e^{-K_{\phi,2}t}W_{2}(\mu, \nu).$$
\end{theorem}

\begin{proof}
By \cite[Theorem~4.4.2]{Wbook2}, it suffices to only consider $\mu=\delta_x$ and $\nu=\delta_y$.
Let $\phi$ be a smooth function in $\mathscr{D}$ and recall that $L=\sum_{i=1}^d(\phi^{-1}V_i)^2+\phi^{-2}Z'$ for the manifold $(M,g')$ and $\{V_i\}_{i=1}^d$ the  $g'$-orthonormal basis of $T_xM$, where $g'=\phi^{-2}g$. Let $X_t$ and $Y_t$ solve the following SDEs respectively:
\begin{align*}
&\vd_I X_t=\sqrt{2}\phi^{-1}(X_t)u_t\,\vd B_t+\phi^{-2}(X_t)Z'(X_t)\,\vd t+N'(X_t)\,\vd l_t,\quad X_0=x;\\
&\vd _I Y_t=\sqrt{2}\phi^{-1}(Y_t)P'_{X_t,Y_t}u_t\,\vd B_t +\phi^{-2}(Y_t)Z'(Y_t)\,\vd t+N'( Y_t)\,\vd \tilde{l_t}, \quad Y_0=y.
\end{align*}
Then, as explained in the proof of Theorem \ref{main-th1}, we have
\begin{align*}%\label{Ito-rho-2}
\vd \rho'(X_t,Y_t)&\leq  \sqrt{2}(\phi^{-1}(X_t)-\phi^{-1}(Y_t))\l<\nabla'\rho'(\newdot,Y_t)(X_t),u_t\,\vd B_t\r>'\notag\\
&\quad-\l(\int_0^{\rho'(X_t,Y_t)}[\phi^{-2}\Ric^Z(\dot{\gamma}(s),\dot{\gamma}(s))+L\log \phi-|\nabla \log \phi|^2](\gamma(s))\,\vd s\r)\vd t.
\end{align*}
Therefore
\begin{align*}
\vd \rho'(X_t, Y_t)^2&=2\rho'(X_t,Y_t)\,\vd \rho'(X_t,Y_t)+(\phi^{-1}(X_t)-\phi^{-1}(Y_t))^2\,\vd t\\
&\leq \vd \tilde{M}_t +2\l[\int_0^{\rho'(X_t,Y_t)}\langle \nabla'\phi^{-1}(\gamma(s)), \dot{\gamma}(s) \rangle' \,\vd s\r]^2\vd t\\
&\quad -2\rho'(X_t,Y_t)\l[\int_0^{\rho'(X_t,Y_t)}[\phi^{-2}\Ric^Z(\dot{\gamma}(s),\dot{\gamma}(s))+L\log \phi-|\nabla \log \phi|^2](\gamma(s))\,\vd s\r]\vd t\\
&\leq \,\vd \tilde{M}_t -2\rho'(X_t,Y_t)\bigg[\int_0^{\rho'(X_t,Y_t)}\lbrace\phi^{-2}(\gamma(s))\Ric^Z(\dot{\gamma}(s),\dot{\gamma}(s))\\
&\quad +L\log \phi(\gamma(s))-2|\nabla\log \phi(\gamma(s))|^2\rbrace\,\vd s\bigg]\,\vd t\\
&\leq \,\vd \tilde{M}_t-2K_{\phi,2}\rho'(X_t,Y_t)^2\,\vd t,
\end{align*}
where
$$\vd \tilde{M}_t=2\sqrt{2}\rho'(X_t,Y_t)(\phi^{-1}(X_t)-\phi^{-1}(Y_t))\l<\nabla'\rho'(\newdot,Y_t)(X_t),u_t\,\vd B_t\r>'.$$
It follows that
\begin{align*}
W_2(\delta_x P_t,\delta_y P_t)^2&\leq \E^{(x,y)}[\rho(X_{t},Y_{t})^2]\leq \|\phi\|_{\infty}^2\E^{(x,y)}[\rho'(X_{t},Y_{t})^2] \\
&\leq \|\phi\|^2_{\infty}\e^{-2K_{\phi,2}t}\rho'(x,y)^2\leq \|\phi\|^2_{\infty}\e^{-2K_{\phi,2} t}\rho(x,y)^2
\end{align*}
which completes the proof.
\end{proof}

We now investigate Talagrand-type inequalities with respect to the uniform distance on the path space $W^{T}:=C([0,T];M)$ of the (reflecting) diffusion process, for a given positive constant $T$. Let $X^{\mu}_t$ be the (reflecting if $\partial M \neq \varnothing$) diffusion process generated by $L$ with initial distribution $\mu\in \mathscr{P}(M)$.
Let $\Pi^{T}_{\mu}$ be the distribution of $$X^{\mu}_{[0,T]}:=\{X^{\mu}_t\colon t\in [0,T]\},$$ which is a probability measure on the (free) path space $W^{T}$. When $\mu=\delta _x$ we denote $\Pi_{\delta_x}^{T}= \Pi_x^{T}$ and $X^{\delta_x}_{[0,T]}=X_{[0,T]}^x$. For any non-negative measurable function $F$ on $W^{T}$ such that $\Pi^{T}_{\mu}(F)=1$, one has
\begin{align}\label{def2}
\mu_F^{T}(\vd x):=\Pi^{T}_{x}(F)\mu (\vd x)\in \mathscr{P}(M).
\end{align}
The the uniform distance on $W^{T}$ is given by
$$\rho_{\infty}(\gamma, \eta):=\sup_{t\in [0,T]}\rho(\gamma_t,\eta_t),\quad\gamma,\eta \in W^{T}.$$
Let $W^{\rho_{\infty}}_2$ be the $L^2$-Wasserstein distance (or
$L^2$-transportation cost) induced by $\rho_{\infty}$. In general,
for any $p\in [1,\infty)$ and two probability measures
$\Pi_1,\Pi_2$ on $W^{T}$,
$$W^{\rho_{\infty}}_{p}(\Pi_1,\Pi_2):=\inf_{\pi \in \mathscr{C}(\Pi_1,\Pi_2)}\l\{\iint_{W^{T}\times W^{T}}\rho_{\infty}(\gamma,\eta)^p\pi(\vd \gamma, \vd \eta)\r\}^{1/p}$$
is the $L^p$-Wasserstein distance (or $L^p$-transportation cost) of
$\Pi_1$ and $\Pi_2$, induced by the uniform norm, where
$\mathscr{C}(\Pi_1, \Pi_2)$ is the set of all couplings for $\Pi_1$
and $\Pi_2$.  Moreover, for
$F\geq 0$ with $\Pi^T_{\mu}(F)=1$, let
$$\mu ^T_{F}(\vd
x)=\Pi^T_{x}(F)\mu(\vd x).$$
The following result improves \cite[Theorems 4.1 and 4.2]{W13} or \cite[Theorems 4.5.3 and 4.5.4]{Wbook2}:

\begin{theorem}\label{thm-tran-ineq}
If there exists $\phi\in \mathscr{D}$ and a constant $K_{\phi}$ satisfying
\begin{align*}
\Ric^Z+ L\log \phi-|\nabla \log \phi|^2\geq K_{\phi}
\end{align*}
then
\begin{enumerate}[\rm(i)]
{ \item  for $F\geq 0, \ \Pi^{T}_{\mu}(F)=1$ and $\mu\in \mathscr{P}(M),$
\small
\begin{align*}
W^{\rho_{\infty}}_2(F\Pi^T_{\mu},\Pi^T_{\mu^{T}_{F}})^2\leq  \frac{2\|\phi\|_{\infty}^2}{K_{\phi}}(\e^{2K_{\phi}^+T}-\e^{2K_{\phi}^-T})\inf_{R>0}\Big\lbrace(1+R^{-1}) \exp \l(8(1+R)\|\nabla\log\phi\|^2_{\infty}\r)\Big\rbrace\ \Pi^{T}_{\mu}(F\log F);
\end{align*}
\normalsize
\item [\rm(i')] for $F\geq 0, \ \Pi^{T}_{\mu}(F)=1$ and $\mu\in \mathscr{P}(M),$
\small
\begin{align*}
W^{\rho_{\infty}}_2(F\Pi^T_{\mu},\Pi^T_{\mu^{T}_{F}})^2\leq  \frac{2\|\phi\|_{\infty}^2}{K_{\phi}}(1-\e^{-2K_{\phi}T})\inf_{R>0}\Big\lbrace(1+R^{-1}) \exp \l(8(1+R)\|\nabla\log\phi\|^2_{\infty}\e^{2K_{\phi}^{+}T}\r)\Big\rbrace\ \Pi^{T}_{\mu}(F\log F);
\end{align*}
\normalsize
}
\item  for any $\mu,\nu\in \mathscr{P}(M)$,
\small
\begin{align*}
W^{\rho_{\infty}}_2(\Pi^T_{\mu},\Pi^T_{\nu})\leq 2\|\phi\|_{\infty}\e^{(K_{\phi}^-+\|\nabla \log \phi\|_{\infty})T}W_{2}(\mu,\nu).
\end{align*}
\normalsize
\end{enumerate}
\end{theorem}

\begin{remark}\label{sec3-rem2}
\text{ }
\begin{enumerate}[(a)]
  \item When $\|\nabla \log \phi\|_{\infty}>0$ and $K_{\phi}>0$ the upper bound in (i) is better than that in (i').
  \item When the boundary is convex we can choose $\phi\equiv 1$. In this case $\nabla \log \phi=0$ and the estimate in (i') is consistent with \cite[Theorem 4.4.2 (2)]{Wbook2} for the convex case.
  \item We note that  \cite[Theorem 4.4.2 (6)]{Wbook2} needs to be corrected as follows:
$$W^{\rho_{\infty}}_2(\Pi^T_{\mu},\Pi^T_{\nu})\leq \e^{K^-T}W_{2}(\mu,\nu),$$ where $K$ is the lower bound of Ricci curvature. It is then consistent with Theorem \ref{thm-tran-ineq} (ii) when $\phi\equiv 1$ and the boundary is convex.
\end{enumerate}
\end{remark}

\begin{proof}[Proof of Theorem \ref{thm-tran-ineq}]\ \smallskip

\noindent(i)\ Simply denote $X^x_{[0,T]}=X_{[0,T]}$. Let $F$ be a positive bounded measurable function on $W^T$ such that $\inf F>0$ and $\Pi^T_{x}(F)=1$. Let $$\vd \mathbb{Q}=F(X_{[0,T]})\,\vd \mathbb{P}.$$ Since $\E F(X_{[0,T]})=\Pi^T_{\mu}(F)=1$, $\mathbb{Q}$ is a probability measure on $\Omega$. Then, we conclude that there exists a unique $\mathscr{F}_t$-predictable process $\beta_t$ on $\mathbb{R}^d$  such that
$$F(X_{[0,T]})=\exp\left({\int_0^T\l<\beta_s, \vd B_s\r>-\frac{1}{2}\int_0^T\|\beta_s\|^2\,\vd s}\right)$$
and
\begin{align}\label{eq13}
\int_0^T\E_{\mathbb{Q}} \|\beta_s\|^2\,\vd s=2\E\left[F(X_{[0,T]})\log F(X_{[0,T]})\right].
\end{align}
Then, by the Girsanov theorem, $\tilde{B}_t:=B_t-\int_0^t\beta_s\,\vd s$, $t\in [0,T]$ is a $d$-dimensional Brownian motion under the probability measure $\mathbb{Q}$.

As explained in the proof of \cite[Theorem 4.5.3]{Wbook2}, it suffices to assume $\mu=\delta_{x}$, $x\in M$. In this case, the desired inequality involves
\begin{align*}
\mu_F^T=\delta_x\ \mbox{ and }\ \Pi^{T}_{\mu}(F\log F)=\Pi^{T}_{x}(F\log F).
\end{align*}
Since the diffusion coefficients are non-constant, it is convenient to adopt the It\^{o} differential $\vd_{I}$ for the Girsanov transformation. So the reflecting  $L$ diffusion process $X_t$ can be constructed by solving the It\^{o} SDE
$$\vd _I X_t=\sqrt{2}\phi^{-1}(X_t)u_t\,\vd B_t+\phi^{-2}(X_t)Z'(X_t)\,\vd t+N'(X_t)\,\vd l_t,\quad\ X_{0}=x,$$
where $B_t$ is the $d$-dimensional Brownian motion with natural filtration $\mathscr{F}_t$. Then
\begin{align}\label{eq1}\vd _I X_t=\sqrt{2}\phi^{-1}(X_t)u_t\,\vd \tilde{B}_t+\{\phi^{-2}(X_t)Z'(X_t)+\sqrt{2}\phi^{-1}(X_t)u_t\beta_t\}\,\vd t+N'(X_t)\,\vd l_t,\quad X_0=x\end{align}
and let $Y_t$ solve
\begin{align}\label{eq2}\vd_I Y_t=\sqrt{2}\phi^{-1}(Y_t)P'_{X_t,Y_t}u_t\,\vd \tilde{B}_t+\phi^{-2}(Y_t)Z'(Y_t)\,\vd t+N'(Y_t)\,\vd \tilde{l}_t,\ \ Y_0=x \end{align}
where $l_t$ and $\tilde{l}_t$ are the local times of $X_t$ and $Y_t$ on $\partial M$, respectively. Moreover, for any bounded measurable function $G$ on $W^T$,
$$\E_{\mathbb{Q}}G(X_{[0,T]}):=\E(FG)(X_{[0,T]})=\Pi_x^T(FG).$$
We conclude that the distribution of $X_{[0,T]}$ under $\mathbb{Q}$ coincides  with $F\Pi^T_x$. Therefore,
\begin{align}\label{eq6}
W^{\rho_{\infty}}_2(F\Pi_x^T,\Pi_x^T)^2&\leq \E_{\mathbb{Q}}\l[\rho_{\infty}(X_{[0,T]},Y_{[0,T]})^2\r]
=\E_{\mathbb{Q}}\Big[\max_{t\in [0,T]}\rho(X_t,Y_t)^2\Big]\notag\\
&\leq \|\phi\|_{\infty}^2\,\E_{\mathbb{Q}}\Big[\max_{t\in [0,T]}\rho'(X_t,Y_t)^2\Big].
\end{align}
Note that due to the convexity of the boundary,
$$\l<N'(x),\nabla'\rho(\newdot,y)(x)\r>'\leq 0,\quad x\in \partial M.$$
From this and equations (\ref{eq1}) and (\ref{eq2}), it follows that
\begin{align*}
\vd \rho'(X_t,Y_t)\leq& \sqrt{2}(\phi^{-1}(X_t)-\phi^{-1}(Y_t))\l<\nabla'\rho'(\newdot,Y_t)(X_t), u_t\,\vd \tilde{B}_t\r>'-K_{\phi}\rho'(X_t,Y_t)\,\vd t+\sqrt{2}\|\beta_t\|\,\vd t.
\end{align*}
Since
$$M_t:=\sqrt{2}\int_0^t\e^{K_{\phi}s}(\phi^{-1}(X_s)-\phi^{-1}(Y_s))\l<\nabla'\rho'(\newdot, Y_s)(X_s),u_s\,\vd \tilde{B}_s\r>'$$
is a $\mathbb{Q}$-martingale, we have
\begin{align*}
\rho'(X_t,Y_t)\leq \e^{-K_{\phi} t}\l(M_t+\sqrt{2}\int_0^t\e^{K_{\phi}s}\|\beta_s\|\,\vd s\r),\quad t\in [0,T].
\end{align*}
So to prove (i), we will estimate the function
$$h_t=\E_{\mathbb{Q}}\max_{s\in [0,t]}\e^{2K_{\phi}s}\rho'(X_s,Y_s)^2.$$
By the Doob inequality, for any $R>0$, we have
\begin{align}\label{Doob}
h_t:&=\E_{\mathbb{Q}}\Big[\max_{s\in [0,t]}\e^{2K_{\phi}s}\rho'(X_s,Y_s)^2\Big]\nonumber\\
&\leq (1+R)\E_{\mathbb{Q}}\max_{s\in [0,t]}M_s^2+2(1+R^{-1})\max_{s\in [0,t]}\E_{\mathbb{Q}}\l[
\l(\int_0^s\e^{K_{\phi}r}\|\beta_r\|\,\vd r\r)^2\r]\nonumber\\
&\leq 4(1+R)\E_{\mathbb{Q}}M_t^2+2(1+R^{-1})\int_0^t\e^{2K_{\phi}s}\,\vd s\int_0^t\E_{\mathbb{Q}}\|\beta_s\|^2\,\vd s\nonumber\\
&\leq  8(1+R)\|\nabla\log\phi\|^2_{\infty}\int_0^th_s\,\vd s
+2(1+R^{-1})\int_0^T\e^{2K_{\phi}s}\,\vd s\,\int_0^T\E_{\mathbb{Q}}\|\beta_s\|^2\,\vd s,\quad t\in [0,T].
\end{align}
 Since $h_0=0$, by using the Gronwall inequality, this inequality further implies
\begin{align}\label{eq7}
 h_T&\leq 2(1+R^{-1}) \exp \l(8(1+R)\|\nabla\log\phi\|^2_{\infty}\r)\int_0^T\e^{2K_{\phi}s}\,\vd s \int_0^T\E_{\mathbb{Q}}\|\beta_s\|^2\,\vd s
\end{align}
By (\ref{eq13}) and \eqref{eq7} we thus have
\begin{align*}
\E_{\mathbb{Q}}\Big[\max_{s\in [0,T]}\rho'(X_s,Y_s)^2\Big]\leq 4(1+R^{-1}) \exp \l(8(1+R)\|\nabla\log\phi\|^2_{\infty}\r)\frac{\e^{2K_{\phi}^+T}-\e^{2K_{\phi}^-T}}{2K_{\phi}}
 \Pi_{x}^T(F\log F).
\end{align*}

\noindent(i')\  \ For this we use the function
$$\tilde{h}_t=\e^{2K_{\phi}t}\E_{\mathbb{Q}}\Big[\max_{s\in [0,t]}\rho'(X_s,Y_s)^2\Big].$$
The inequality \eqref{Doob} should then be modified as follows:
\begin{align*}
\tilde{h}_t&:=\e^{2K_{\phi}t}\E_{\mathbb{Q}}\Big[\max_{s\in [0,t]}\rho'(X_s,Y_s)^2\Big]\nonumber\\
&\,\,\leq \e^{2K_{\phi}t}\l((1+R)\E_{\mathbb{Q}}\Big[\max_{s\in [0,t]}\e ^{-2K_{\phi}s}M_s^2\Big]+2(1+R^{-1})\max_{s\in [0,t]}\E_{\mathbb{Q}}\l[\l(\int_0^s\e^{-K_{\phi}(s-r)}\|\beta_r\|\,\vd r\r)^2\r]\r)\nonumber\\
&\,\,\leq4(1+R)\e^{2K_{\phi}^+t}\E_{\mathbb{Q}}M_t^2+2(1+R^{-1})\int_0^t\e^{2K_{\phi}r}\,\vd r\int_0^t\E_{\mathbb{Q}}\|\beta_s\|^2\,\vd s\nonumber\\
&\,\,\leq 8(1+R)\|\nabla\log\phi\|^2_{\infty}\e^{2K_{\phi}^+T}\int_0^t\tilde{h}_s\,\vd s
+2(1+R^{-1})\int_0^T\e^{2K_{\phi}r}\,\vd r\int_0^T\E_{\mathbb{Q}}\|\beta_s\|^2\,\vd s,\quad\ t\in [0,T].
\end{align*}
Since $\tilde{h}_0=0$, this inequality implies
\begin{align*}
\tilde{h}_T&\leq 2(1+R^{-1}) \exp \l(8(1+R)\|\nabla\log\phi\|^2_{\infty}\e^{2K_{\phi}^+T}\r)\int_0^T\e^{2K_{\phi}s}\,\vd s \int_0^T\E_{\mathbb{Q}}\|\beta_s\|^2\,\vd s.
\end{align*}
We then conclude that
\begin{align*}
\E_{\mathbb{Q}}\max_{s\in [0,T]}\rho'(X_s,Y_s)^2\leq 4(1+R^{-1}) \exp \l(8(1+R)\|\nabla\log\phi\|^2_{\infty}\e^{2K_{\phi}^+T}\r)\frac{1-\e^{-2K_{\phi}T}}{2K_{\phi}}
 \Pi_{x}^T(F\log F).
\end{align*}

{\setlength{\parindent}{0cm}
(ii) Without loss of generality, we consider $\mu=\delta_x$, and $\nu=\delta_y$. Let $X_t$ and $Y_t$ solve the following SDEs, respectively:
}
\begin{align*}
&\vd_I X_t=\sqrt{2}\phi^{-1}(X_t)u_t\,\vd B_t+\phi^{-2}(X_t)Z'(X_t)\,\vd t+N'(X_t)\,\vd l_t,\quad X_0=x;\\
&\vd _I Y_t=\sqrt{2}\phi^{-1}(Y_t)P'_{X_t,Y_t}u_t\,\vd B_t +\phi^{-2}(Y_t)Z'(Y_t)\,\vd t+N'( Y_t)\,\vd \tilde{l_t}, \quad Y_0=y.
\end{align*}
Then, as explained in the proof of Theorem \ref{main-th1}, we have
\begin{align}\label{Ito-rho-1}
\vd \rho'(X_t,Y_t)&\leq  \sqrt{2}(\phi^{-1}(X_t)-\phi^{-1}(Y_t))\l<\nabla'\rho'(\newdot,Y_t)(X_t),u_t\,\vd B_t\r>'\notag\\
&\quad-\l(\int_0^{\rho'(X_t,Y_t)}\big(\phi^{-2}\Ric^Z(\dot{\gamma}(s),\dot{\gamma}(s))+L\log \phi-|\nabla \log \phi|^2\big)(\gamma(s)\,\vd s\r)\vd t.
\end{align}
Therefore,
\begin{align}\label{Ieq2}
\rho'(X_t,Y_t)\leq \e^{-K_{\phi}t}(\hat{M}_t+\rho'(x,y)),\quad t\geq 0
\end{align}
for
$$\hat{M}_t:=\sqrt{2}\int_{0}^t\e^{K_{\phi}s}(\phi^{-1}(X_s)-\phi^{-1}(Y_s))\l<\nabla'\rho_{s}(\newdot, Y_s)(X_s),u_s\,\vd B_s\r>'.$$
Again using the It\^{o} formula, we have
\begin{align*}
\vd \rho'(X_t,Y_t)^2\leq \vd \tilde{M}_t-2(K_{\phi}-\|\nabla \log \phi\|^2_{\infty})\rho'(X_t,Y_t)^2\,\vd t
\end{align*}
where
$$\vd \tilde{M}_t=2\rho'(X_t,Y_t)(\phi^{-1}(X_t)-\phi^{-1}(Y_t))\l<\nabla'\rho'(\newdot,Y_t)(X_t),u_t\,\vd B_t\r>'$$
which implies
\begin{align*}
\E \rho'(X_{t},Y_{t})^2\leq \e^{-2(K_{\phi}-\|\nabla \log \phi\|_{\infty}^2)t}\rho'(x,y)^2.
\end{align*}
Combining this with \eqref{Ieq2} we arrive at
\begin{align*}
W_2^{\rho_{\infty}}(\Pi^T_{x},\Pi^T_{y})^2&\leq \|\phi\|^2_{\infty} \,\E\max_{t\in [0,T]}\rho'(X_t,Y_t)^2\\
&\leq \|\phi\|^2_{\infty} \e^{2K_{\phi}^-T}\E\max_{t\in [0,T]}(\hat{M}_t+\rho'(x,y))^2\\
&\leq 4 \|\phi\|^2_{\infty} \e^{2K_{\phi}^-T}\E[\hat{M}_{T}+\rho'(x,y)]^2=4\|\phi\|^2_{\infty}\,\e^{2K_{\phi}^-T}\l(\E[\hat{M}_{T}^2]+\rho'(x,y)^2\r)\\
&\leq 4 \|\phi\|^2_{\infty} \e^{2K_{\phi}^-T}\l(2\int_{0}^{T}\e^{2K_{\phi}t}\|\nabla \log \phi\|^2_{\infty}\,\E\rho'(X_t,Y_t)^2\,\vd t+\rho'(x,y)^2\r)\\
&\leq 4 \|\phi\|^2_{\infty} \e^{2(K_{\phi}^-+\|\nabla \log\phi\|^2_{\infty})T}\rho'(x,y)^2\\
&\leq 4 \|\phi\|^2_{\infty} \e^{2(K_{\phi}^-+\|\nabla \log\phi\|^2_{\infty})T}\rho(x,y)^2
\end{align*}
where the second inequality is due to the Doob inequality. This implies the desired inequality for
$\mu=\delta_x$ and $\nu=\delta_y$.
\end{proof}

\begin{corollary}\label{sec3-cor-2}
If there exists $\phi\in \mathscr{D}$ and a constant $K_{\phi}$ satisfying
\begin{align*}
\Ric^Z+ L\log \phi-|\nabla \log \phi|^2\geq K_{\phi}
\end{align*}
then
\begin{enumerate}[\rm(i)]
\item  for $F\geq 0$,  $\Pi^{T}_{\mu}(F)=1$ and $\mu\in \mathscr{P}(M),$
\small
\begin{align*}
W^{\rho_{\infty}}_2(F\Pi^T_{\mu},\Pi^T_{\mu^{T}_{F}})^2\leq  \frac{2\|\phi\|_{\infty}^2}{K_{\phi}}\l(\e^{2K_{\phi}^+T}-\e^{2K_{\phi}^-T}\r) \exp \l(8\|\nabla\log\phi\|^2_{\infty}+4\sqrt{2}\|\nabla\log\phi\|_{\infty}\r)\ \Pi^{T}_{\mu}(F\log F);
\end{align*}
\normalsize
\item [\rm(i')] for $F\geq 0$, $\Pi^{T}_{\mu}(F)=1$ and $\mu\in \mathscr{P}(M),$
\small
\begin{align*}
W^{\rho_{\infty}}_2(F\Pi^T_{\mu},\Pi^T_{\mu^{T}_{F}})^2\leq  \frac{2\|\phi\|_{\infty}^2}{K_{\phi}}\l(1-\e^{-2K_{\phi}T}\r)\exp \l(8\|\nabla\log\phi\|^2_{\infty}\e^{2K_{\phi}^{+}T}+4\sqrt{2}\|\nabla\log\phi\|_{\infty}\e^{K_{\phi}^{+}T}\r)\ \Pi^{T}_{\mu}(F\log F).
\end{align*}
\normalsize
\end{enumerate}
\end{corollary}
\begin{proof}
It is easily observed that
\begin{align*}
(1+R^{-1})\exp\l(8(1+R)\|\nabla\log\phi\|^2_{\infty}\r)\leq
\exp \l(R^{-1}+8(1+R)\|\nabla\log\phi\|^2_{\infty}\r).
\end{align*}
Taking the infimum about $R$  on the right side above, we arrive at
\begin{align*}
 \exp \l(R^{-1}+8(1+R)\|\nabla\log\phi\|^2_{\infty}\r)\geq   \exp \l(8\|\nabla\log\phi\|^2_{\infty}+4\sqrt{2}\|\nabla \log \phi\|_{\infty}\r)
\end{align*}
which allows to prove (i). The inequality (i') can be checked in the same way.
\end{proof}

\section{New construction of function $\log \phi$}

In this section, we give a new construction of a function $\phi$ which satisfies the conditions of the previous section. To do so, we let $\rho_{\partial}$ be the Riemannian distance to the boundary $\partial M$ and use a comparison theorem for $\Delta\rho_{\partial}$ near the boundary, essentially due to \cite{Kasue}. Note that, by using local charts, it is clear that $\rho_{\partial}$ is smooth in a neighbourhood of $\partial M$. We call
$$i_{\partial}:=\sup\{r>0: \rho_{\partial} \mbox{ is smooth on } \{\rho_{\partial}<r\}\}$$
the injectivity radius of $\partial M$. Obviously, $i_{\partial}>0$ if $M$ is compact, but it could be zero in the non-compact case ($\sup \varnothing =0$ by convention). As \cite[Theorem~1.2.3]{Wbook2} we have:

\begin{lemma}\label{lem1-sec3}
Let $\theta, k$ be constants such that $\II \leq \theta$ and $\Sect  \leq k$. Let
\begin{align}\label{fun-h}
h(t):=\left\{
       \begin{array}{ll}
         \cos{\sqrt{k}t}-\frac{\theta}{\sqrt{k}}\sin {\sqrt{k}t}, & k\geq 0, \\
         \cosh{\sqrt{-k}t}-\frac{\theta}{\sqrt{-k}}\sinh{\sqrt{-k}t}, & k<0
       \end{array}
     \right.
\end{align}
for $t\geq 0$. Let $h^{-1}(0)$ be the first zero of $h$ (with $h^{-1}(0):=\infty$ if $h(t)>0$ for all $t\geq0$).
Then for any $x\in \mathring{M}$ such that $\rho_{\partial}(x)\leq i_{\partial}\wedge h^{-1}(0)$ we have
\begin{align}\label{compare-thm}
\Delta \rho_{\partial}(x)\geq (d-1)\frac{h'}{h}(\rho_{\partial}(x)).
\end{align}
\end{lemma}

Note that if $k$ is positive then
$$h^{-1}(0) = \frac{1}{\sqrt{k}}\arcsin\l(\sqrt{\frac{k}{k+\theta^2}}\r).$$
We now work under the following assumption:
\smallskip

\noindent{\bf Assumption (A) }
There exist non-negative constants $\sigma$ and $\theta$ such that $-\sigma\leq \II \leq \theta$ and a positive constant $r_0$ such that on $\partial_{r_0}M:=\{x\in M: \rho_{\partial}(x)\leq r_0\}$ the function $\rho_{\partial}$ is smooth, the norm of $Z$ is bounded and $\Sect \leq k$ for some positive constant $k$.
\medskip

Using this assumption, F.-Y. Wang constructed a function $\phi$ satisfying $\phi\in \mathscr{D}$
   (see \cite[p.1436]{W07} or \cite[Theorem 3.2.9]{Wbook2}).
Following his construction, we define
\begin{align*}
\log \phi(x) =\frac{\sigma}{\alpha}\int_0^{\rho_{\partial}(x)}[h(s)-h(r_1)]^{1-d}\vd s \int_{s\wedge r_1}^{r_1}[h(u)-h(r_1)]^{d-1}\vd u,
\end{align*}
where $r_1:=r_0\wedge h^{-1}(0)$ and $$ \alpha:=(1-h(r_1))^{1-d}\int_0^{r_1}[h(s)-h(r_1)]^{d-1}\vd s.$$ Then from the proof of \cite[Theorem 1.1]{Wa05}, we know:

\begin{theorem}
Suppose that Assumption ${\bf (A)}$ holds and $\Ric^Z\geq K$.
Define
\begin{align*}
K_p&=K-\sigma\left(\delta_{r_1}(Z)+\frac{d}{r_1}\right)-p\sigma^2,
\end{align*}
where \begin{equation}\label{Eq:r_1}
\delta_{r_1}(Z):=\sup \big\lbrace |Z(x)|\colon
x\in \partial_{r_1}M \big\rbrace.
\end{equation}
Then all results in Section 2  hold by replacing
$$K_{\phi},\ K_{\phi,p},\  \|\phi\|_{\infty}\  \mbox{ and }\ \|\nabla \log \phi\|_{\infty}$$
with
$$K_1, \ {K}_{p}, \ \e^{\frac{1}{2}\sigma d r_1}\  \mbox{ and }\ \sigma$$
respectively.
\end{theorem}

    In the following we
   give a new construction of function $\phi$ by using the modifier
   \begin{align*}
\ell(r)=\begin{cases}
\e^{-2}-\e^{-{2}(1-2r)^{-1}}, & 0\leq r<\frac12, \\
\e^{-2}, & r\geq \frac12.
\end{cases}
\end{align*}

\begin{proposition}\label{sec3-prop}
Suppose that Assumption ${\bf (A)}$ holds. Let
$$H(r):=\frac{\sqrt{k+\theta^2}}{k}\cos\left(\arcsin\left(\sqrt{\frac{k}{k+\theta^2}}\right)
-\sqrt{k}\,(r\wedge r_1)\right)-\frac{\theta}{k}.$$
Then the function
\begin{align}\label{log-phi}
\log \phi (x):=\frac{1}{2}\sigma \e^2\ell\l(\frac{H(\rho_{\partial}(x))}{2H(r_1)}\r)H(r_1)
\end{align}
satisfies
$$N\log \phi|_{\partial M}=\sigma\geq -\II .$$
Moreover,
\begin{align*}
\|\phi\|_{\infty}\leq \e^{\frac{1}{2}\sigma H(r_1)}, \quad |\nabla \log \phi|\leq \sigma
\end{align*}
and
\begin{align*}
L \log \phi(x)\geq -\sigma\left(d\sqrt{\theta^2+k}+\delta_{r_1}(Z)+\frac{5}{2 H(r_1)}\right).
\end{align*}

\end{proposition}

\begin{proof}[Proof of Proposition \ref{sec3-prop}]
First it is easy to see that the modifier $\ell$ satisfies
$\ell\leq \e^{-2}$.
Differentiating $\ell$ we obtain
\begin{align*}
\ell'(r)=\begin{cases}
         \big(\frac12-r\big)^{-2}\e^{-\big(\frac12-r\big)^{-1}}, & 0\leq r<\frac12;\\
             \ 0, & r\geq\frac12
           \end{cases}
 \end{align*}
and
\begin{align*}
\ell''(r)=\begin{cases}
              -{2r}\,\big(\frac12-r\big)^{-4}\e^{-\big(\frac12-r\big)^{-1}}, & 0\leq r<\frac12; \\
            \quad  0, & r\geq \frac12.
            \end{cases}
\end{align*}
As $\ell''<0$ on $[0,1/2)$, the function $\ell'$ is at its maximal point when $r=0$, which implies $0\leq  \ell'\leq 4\e^{-2}$. Using the same method, when $r={\sqrt3}/6$ the function $\ell''$ reaches the minimal value, which implies
$$\ell''\geq -3^{-1/2}(3+\sqrt{3})^4\e^{-(3+\sqrt{3})}> -20\e^{-2}.$$
Using these results, we have
\begin{align*}
N\log \phi|_{\partial M}=\frac{1}{4}\e^2\sigma \ell'(0)N\rho_{\partial}=\sigma,
\end{align*}
and
\begin{align*}
|\nabla \log \phi|=\frac{1}{4}\e^{2}\sigma\ell'\l(\frac{H(\rho_{\partial})}{2H(r_0)}\r)H'(\rho_{\partial})\leq \sigma.
\end{align*}
Moreover, by Lemma \ref{lem1-sec3}, we have
\begin{align*}
L\log \phi&=\frac{1}{4}\e^{2}\sigma\left(\ell'\l(\frac{H(\rho_{\partial})}{2H(r_0)}\r)h(\rho_{\partial})L\rho_{\partial} +\ell''\l(\frac{H(\rho_{\partial})}{2H(r_0)}\r)\frac{h(\rho_\partial)^2}{2H(r_0)}+\ell'\l(\frac{H(\rho_{\partial})}{2H(r_0)}\r)h'(\rho_{\partial}(x))\right)\\
&\geq \frac{1}{4}\e^{2}\sigma\left(\ell'\l(\frac{H(\rho_{\partial})}{2H(r_0)}\r)\Big(dh'(\rho_{\partial})-\sup_{\partial_{r_0}M}|Z|\Big) +\frac{h(\rho_\partial)^2}{2H(r_0)}\ell''\l(\frac{H(\rho_{\partial})}{2H(r_0)}\r)\right),
\end{align*}
where $h$ is defined as in \eqref{fun-h} for $k\geq 0$. It is easy to calculate that
\begin{align*}
h'(r)&\geq -\sqrt{\theta^2+k}.
\end{align*}
Combining this with properties of $\ell$, we conclude that
\begin{align*}
L\log \phi\geq&- \sigma\left(d\sqrt{\theta^2+k}+\sup\left\lbrace |Z|(x)\colon x \in \partial_{r_0\wedge h^{-1}(0)}M \right\rbrace+\frac{5}{2 H(r_0)}\right)
\end{align*}
which completes the proof.
\end{proof}

\begin{corollary}\label{sec3-cor-1}
Suppose that Assumption ${\bf (A)}$ holds and $\Ric^Z\geq K$.
Define
\begin{align*}
\tilde{K}&=K-\sigma\left(d\sqrt{\theta^2+k}+\delta_{r_1}(Z)+\frac{5}{2 H(r_1)}\right),
\end{align*}
and $\tilde{K}_{p}=\tilde{K}-p\sigma^2$ with  $\delta_{r_1}(Z)$ as defined in \eqref{Eq:r_1}.
Then all results in Section 2  hold by replacing
$$K_{\phi,p},\  \|\phi\|_{\infty}\  \mbox{ and }\ \|\nabla \log \phi\|_{\infty}$$
with
$$\tilde{K}_{p}, \ \e^{\frac{1}{2}\sigma H(r_1)}\  \mbox{ and }\ \sigma$$
respectively.
\end{corollary}

\proof[Acknowledgements]This work has been supported by Fonds National
de la Recherche Luxembourg (Open project O14/7628746 GEOMREV). The
first named author acknowledges support by NSFC (Grant No.~11501508)
and Zhejiang Provincial Natural Science Foundation of China (Grant
No. LQ16A010009).

\providecommand{\bysame}{\leavevmode\hbox to3em{\hrulefill}\thinspace}
\providecommand{\MR}{\relax\ifhmode\unskip\space\fi MR }
\providecommand{\MRhref}[2]{%
  \href{http://www.ams.org/mathscinet-getitem?mr=#1}{#2}
}
\providecommand{\href}[2]{#2}

\end{document}